\documentclass[review]{elsarticle}

\usepackage{multirow}
\usepackage{subfigure}
\usepackage{amsmath,xfrac,cases,mathrsfs} 
\usepackage{amssymb}
\usepackage{enumerate}
\usepackage{lineno,hyperref}
\usepackage{braket}

\newtheorem{thm}{Theorem}
\newtheorem{lem}[thm]{Lemma}
\newdefinition{rmk}{Remark}
\newdefinition{prop}[thm]{Proposition}
\newproof{pf}{Proof}
\newproof{pot}{Proof of Theorem \ref{thm2}}

\journal{!}

\begin{document}
	
\begin{frontmatter}
	
\title{Nadaraya-Watson estimator for reflected stochastic processes driven by Brownian motions}

\author[rvt]{Han Yuecai}
\ead{hanyc@jlu.edu.cn}
		
\author[rvt]{Zhang Dingwen\corref{cor1}}
\ead{zhangdw20@mails.jlu.edu.cn}

\cortext[cor1]{Corresponding author}

\address[rvt]{School of Mathematics, Jilin University, Changchun, 130000, China}
		
\begin{abstract}
We study the Nadaraya-Watson (N-W) estimator for the drift function of two-sided reflected stochastic processes. We propose a discrete-type N-W estimator and a continuous-type N-W estimator based on the discretely observed processes and continuously observed processes respectively. Under some regular conditions, we obtain the consistency and establish the asymptotic distributions for the two estimators. Furthermore, we briefly remark that our method can be applied to the one-sided reflected stochastic processes spontaneously.  Numerical studies show that the proposed estimators is adequate for practical use.
\end{abstract}
		
\begin{keyword}
Reflected stochastic processes \sep Discrete observations \sep Continuous observations \sep Nadaraya-Watson estimator 
\end{keyword}
		
\end{frontmatter}
	
\section{Introduction}\label{sec1}
Let $(\Omega, \mathcal{F}, \mathbb{P}, \{\mathcal{F}\}_{t\geq0})$ be a filtered probability space and $W=\{W_{t}\}_{t\geq0}$ is a one-dimensional standard Brownian motion adapted to $\{\mathcal{F}_{t}\}_{t\geq0}$. The stochastic process $X=\{X_{t}\}_{t\geq0}$  is defined as the unique strong solution to the following reflected stochastic differential equation (SDE) with two-sided barriers $l$ and $u$
\begin{equation}\label{eqmodelT}
\left\{
\begin{aligned}
&\mathrm{d}X_{t}=b(X_{t})\mathrm{d}t+\sigma\mathrm{d}W_{t}+\mathrm{d}L_{t}-\mathrm{d}R_{t},\\
&X_{0}=x\in[l,u],
\end{aligned}
\right.
\end{equation}
where $\sigma>0$, $0\leq l<u<\infty$ are given numbers, and $b:[l,u]\rightarrow\mathbb{R}$ is an unknown measurable function. The regular conditions on the drift function $b$ will be provided on Section \ref{sec2}. The solution of a reflected SDE behaves like the solution of a standard SDE in the interior of its domain $(l,u)$. The processes $L=\{L_{t}\}_{t\geq0}$ and $R=\{R_{t}\}_{t\geq0}$ are the nimimal continuous increasing processes such that $X_{t}\in[l,u]$ for all $t\geq0$. Moreover, the two processes $L$ and $R$ subject to $L_{0}=R_{0}=0$ increase only when $X$ hits its boundary $l$ and $u$, and
\begin{equation*}
\int_{0}^{\infty}I_{\{X_{t}>l\}}\mathrm{d}L_{t}=0, \quad \int_{0}^{\infty}I_{\{X_{t}<u\}}\mathrm{d}R_{t}=0,
\end{equation*} 
where $I(\cdot)$ is the indicator function. For more about reflected stochastic processes, one can refer to \cite{Harrison2013}.
Assume that the discretely observed process $X$ is observed at regularly sapced time points $\{t_{k}=k\Delta, k=0,1, \cdots, n\}$. We shall discuss the two types of N-W estimator of drift function $b$ based on the discrete observations $\{X_{t_{k}},L_{t_{k}},R_{t_{k}}, k=0,1,\cdots,n\}$ and continuous observations $\{X_{t},L_{t},R_{t}\}_{\{t\geq0\}}$ respectively.

Reflected SDEs have been widely used in many fields such as the queueing system \citep{Ward2003a,Ward2003b,Ward2005}, financial engineering \citep{Bo2010,Bo2011a,Han2016} and mathematical biology \citep{Ricciardi1987}. The reflected barrier is usually bigger than $0$ due to the physical restriction of the state processes which take non-negative values. One can refer to \cite{Harrison1985,Whitt2002} for more details on reflected SDEs and their broad applications. 

The parameter estimation problem in reflected SDEs has gained much attention in recent years. A drift parameter estimator for a reflected fractional Brownian motion is proposed in \cite{Hu2013}. For reflected Ornstein–Uhlenbeck (OU) processes, on the one hand,
a maximum likelihood estimator (MLE) for the drift parameter based on the continuously observed processes is proposed in \cite{Bo2011b}. A sequential MLE for the drift parameter based on the continuous observations throughout a random time interval is proposed in \cite{Lee2012}. The heart of MLE is the Girsanov theorem. On the other hand, an ergodic type estimator for drift parameters based on discretely observed reflected OU processes is proposed in \cite{Hu2015}. Subsequently, an ergodic type estimator for all parameters (drift and diffusion parameters) is proposed in \cite{Hu2021}. However, there is only limited literature on the nonparametric estimation for the drift function of reflected SDEs. 

N-W estimator has been widely used on nonparametric estimation of stochastic processes. There are several types of N-W estimator For stochastic processes driven by fractional Brownian motion \cite{Prakasa2011,Saussereau2014,Fabienne2019}. For stochastic processes driven by stable L\'evy motions, a nonparametric N-W estimator is proposed \cite{Long2013}. For more about N-W estimator, one can refer to \cite{N1964,W(1964)}.

The remainder of this paper is organized as follows. In section \ref{sec2}, we describe some notations and assumptions related to our context. Under the assumptions, we propose the discrete-type and continuous-type N-W estimator, obtain the consistency and establish the asymptotic distributions of the two estimators. We also discuss the extension of main results to one-sided reflected processes. The proofs of main results are showed in Section \ref{sec3}. In section \ref{sec4}, we present some examples and give
their numerical results. Section \ref{sec5} concludes with some discussion and remarks on the further work.

\section{Main Results}\label{sec2}
\subsection{Notation and Assumptions}
In this subsection, we present the necessary conditions and the construction of the estimator. Throughout the paper, we shall use notation ``$\stackrel{P}{\longrightarrow}$" to denote ``convergence in probability" and notation ``$\sim$" to denote ``convergence in distribution”. Let $\|f\|_{\infty}=\sup_{\cdot\in\mathbb{R}}f(\cdot)$.

Now we introduce the following set of assumptions.
\begin{enumerate}[i.]
	\item\label{asp1} The drift function $b(\cdot)$ satisfies a global Lipschitz condition, i.e., there exists a positive constant $M>0$ such that
	\begin{equation*}
	|b(x_{1})-b(x_{2})|\leq M|x_{1}-x_{2}|, \quad x_{1}, x_{2}\in [l,u].
	\end{equation*}
    \item\label{asp2}$b(\cdot)$ is not identically $0$ on $[l,u]$.
    \item\label{asp3} $K_{h}(\cdot)=K(\cdot/h)/h$, where $h$ is the bandwidth, and the kernel function $K(\cdot)$ is a symmetric probability density with support $[—1, 1]$, mean 0, and bounded first derivative.
    \item\label{asp4} For discretely observed processes, there exists an $\epsilon$ small enough, such that $\Delta\rightarrow0$, $h\rightarrow0$, $n\Delta\rightarrow\infty$, and $\Delta^{\frac{1}{2}-\epsilon}h^{-1}\rightarrow0$, as $n\rightarrow\infty$.
    \item\label{aspl} For continuously observed processes, $h\rightarrow0$, as $T\rightarrow\infty$.
\end{enumerate}

Under condition \ref{asp1}, there exists a unique solution of Eq. (\ref{eqmodelT}) by an extension results of \cite{Lions1984}.

The N-W estimator based on discretely observed processes is given by
\begin{equation*}
\hat{b}_{n}(x)=\frac{\sum_{k=0}^{n-1}K_{h}(X_{t_{k}}-x)(\delta X_{t_{k}}-\delta L_{t_{k}}+\delta R_{t_{k}})}{\Delta\sum_{k=0}^{n-1}K_{h}(X_{t_{k}}-x)},
\end{equation*}
while based on continuously observed processes is given by
\begin{equation*}
\tilde{b}_{T}(x)=\frac{\int_{0}^{T}K_{h}(X_{t}-x)(\mathrm{d}X_{t}-\mathrm{d}L_{t}+\mathrm{d}R_{t})}{\int_{0}^{T}K_{h}(X_{t}-x)\mathrm{d}t}.
\end{equation*}

To this end, we give the unique invariant and the ergodic property for $X$.
\begin{lem}\label{lemE}
Under conditions \ref{asp1}-\ref{aspl}, for any function $f\in L_{1}([l,u])$, we have,
\begin{enumerate}[a.]
	\item $X$ has a unique invariant probability measure,
	\begin{equation*}
		\pi_{l,u}(x)=\frac{e^{-\frac{2}{\sigma^{2}} \int_{l}^{x} b(y) \mathrm{d} y}}{\int_{l}^{u} e^{-\frac{2}{\sigma^{2}} \int_{l}^{x} b(y) \mathrm{d} y} \mathrm{d} x}.
	\end{equation*}
    \item The continuously observed processes $\{X_{t}\}_{\{t\geq0\}}$ is ergodic,
    \begin{equation*}
    	\lim_{T\rightarrow\infty}\frac{1}{T}\int_{0}^{T}f(X_{t})\mathrm{d}t=\mathbb{E}[f(X_{\infty})]=\int_{l}^{u}f(x)\pi_{l,u}(x)\mathrm{d}x.
    \end{equation*}
	\item The discretely observed processes $\{X_{t_{k}}, k=1,2,\cdots,n\}$ is ergodic, 
	\begin{equation*}
	\lim_{n\rightarrow\infty}\frac{1}{n}\sum_{k=0}^{n-1}f(X_{t_{k}})=\mathbb{E}[f(X_{\infty})]=\int_{l}^{u}f(x)\pi_{l,u}(x)\mathrm{d}x.
	\end{equation*}
\end{enumerate}
\end{lem}

\subsection{Asymptotic bahavior of the N-W estimator}
In this subsection, we propose the consistency of the two types of N-W estimators and establish their rates of convergence and asymptotic distributions. 

We define
\begin{equation}\label{eqFr}
	F(r)=\int_{l}^{u}K_{h}(y-x)\pi_{l,u}(y)\mathrm{d}y=\int_{-1}^{1}K(r)\pi_{l,u}(x+rh)\mathrm{d}r,
\end{equation}
and
\begin{equation}\label{eqSga}
	\Sigma(r)=\sigma^{2}F^{-1}(r).
\end{equation}

The consistency of discrete-type N-W estimator $\hat{b}_{n}$ is given as follows.
\begin{thm}\label{thmConsD}
Under conditions \ref{asp1}-\ref{aspl}, we have 
\begin{equation*}
\hat{b}_{n}(x)\stackrel{P}{\longrightarrow} b(x),
\end{equation*}
as $n\rightarrow\infty$.
\end{thm}

The consistency of continuous-type N-W estimator $\tilde{b}_{T}$ is given as follows.
\begin{thm}\label{thmConsC}
Under conditions \ref{asp1}-\ref{aspl}, we have
\begin{equation*}
\tilde{b}_{T}(x)\stackrel{P}{\longrightarrow}b(x),
\end{equation*}
as $T\rightarrow\infty$.
\end{thm}

To study the asymptotic distribution and rate of convergence of the two types N-W estimator, we impose some new conditions as follow.
\begin{enumerate}[1.]
	\item \label{aspamd} For discretely observed processes, there exists an $\epsilon$ small enough, $nh\Delta\rightarrow\infty$, $nh^{3}\Delta\rightarrow0$ and $nh^{-1}\Delta^{2-\epsilon}\rightarrow0$, as $n\rightarrow\infty$.
	\item \label{aspamc} For continuously obseved processes, $Th\rightarrow\infty$ and $Th^{3}\rightarrow0$, as $T\rightarrow\infty$.
\end{enumerate}

The asymptotic distribution and rate of convergence of the discrete-type N-W estimator for the drift function is given as follows.
\begin{thm}\label{thmAD}
Under conditions \ref{asp1}-\ref{asp4} and \ref{aspamd}, we have
\begin{equation*}
\sqrt{nh\Delta}\big(\hat{b}_{n}(x)-b(x)\big)\sim\mathcal{N}(0,\Sigma(x)),
\end{equation*}
as $n\rightarrow\infty$.
\end{thm}

The asymptotic distribution and rate of convergence of the continuous-type N-W estimator for the drift function is given as follows.
\begin{thm}\label{thmAC}
	Under conditions \ref{asp1}-\ref{asp3}, \ref{aspl} and \ref{aspamc}, we have
	\begin{equation*}
		\sqrt{Th}\big(\tilde{b}_{T}(x)-b(x)\big)\sim\mathcal{N}(0,\Sigma(x)),
	\end{equation*}
	as $n\rightarrow\infty$.
\end{thm}

\subsection{Results of an expansion}
In this subsection, we discuss the extension results of our main results to the reflected SDE with only one-sided barrier. Note that it is almost the same for only a lower reflecting barrier or a upper reflecting barrier. We consider the following reflected SDE with a lower reflecting barrier $l$
\begin{equation}\label{eqmodelo}
	\left\{
	\begin{aligned}
		&\mathrm{d}X_{t}=b(X_{t})\mathrm{d}t+\sigma\mathrm{d}W_{t}+\mathrm{d}L_{t},\\
		&X_{0}=x\in[l,\infty).
	\end{aligned}
	\right.
\end{equation}
Hence, the discrete-type N-W estimator is given by
\begin{equation*}
	\hat{b}_{n}(x)=\frac{\sum_{k=0}^{n-1}K_{h}(X_{t_{k}}-x)(\delta X_{t_{k}}-\delta L_{t_{k}})}{\Delta\sum_{k=0}^{n-1}K_{h}(X_{t_{k}}-x)},
\end{equation*}
while the continuous-type N-W estimator is given by
\begin{equation*}
	\tilde{b}_{T}(x)=\frac{\int_{0}^{T}K_{h}(X_{t}-x)(\mathrm{d}X_{t}-\mathrm{d}L_{t})}{\int_{0}^{T}K_{h}(X_{t}-x)\mathrm{d}t}.
\end{equation*}

\begin{rmk}
Our method can be applied to the reflected processes with only one-sided barrier. The unique invariant density of $X$ is given by
\begin{equation*}
\pi_{l}(x)=\frac{e^{-\frac{2}{\sigma^{2}} \int_{l}^{x} b(y) \mathrm{d} y}}{\int_{l}^{\infty} e^{-\frac{2}{\sigma^{2}} \int_{l}^{x} b(y) \mathrm{d} y} \mathrm{d} x}.
\end{equation*}
With the unique invariant density, we could do similar proofs as the proofs of Theorem \ref{thmConsD}, \ref{thmConsC}, \ref{thmAD} and \ref{thmAC} to establish their consistency and asymptotic distribution. We omit the details here.
\end{rmk}

\section{Proofs of the Main Results}\label{sec3}
In this section, we give the proofs of the main results. We first propose the useful alternative expressions for the two estimators. Note that
\begin{equation*}
\delta X_{t_{k}}=\int_{t_{k}}^{t_{k+1}}b(X_{t})\mathrm{d}t+\sigma\delta W_{t_{k}}+\delta L_{t_{k}}-\delta R_{t_{k}}.
\end{equation*}
Then the discrete-type N-W estimator $\hat{b}_{n}(x)$ turns to
\begin{equation*}
\begin{aligned}
\hat{b}_{n}(x)=&\frac{\frac{1}{n\Delta}\sum_{k=0}^{n-1}K_{h}(X_{t_{k}}-x)\int_{t_{k}}^{t_{k+1}}b(X_{t})\mathrm{d}t}{\frac{1}{n}\sum_{k=0}^{n-1}K_{h}(X_{t_{k}}-x)}\\
&+\frac{\frac{\sigma}{n\Delta}\sum_{k=0}^{n-1}K_{h}(X_{t_{k}}-x)\delta W_{t_{k}}}{\frac{1}{n}\sum_{k=0}^{n-1}K_{h}(X_{t_{k}}-x)}\\
=&\frac{G_{n,1}(x)+G_{n,2}(x)}{F_{n}(x)}.
\end{aligned}
\end{equation*}
And also a useful alternative expression for the continuous-type N-W estimator $\tilde{b}(x)$ 
\begin{equation*}
\begin{aligned}
\tilde{b}_{T}(x)=&\frac{\frac{1}{T}\int_{0}^{T}K_{h}(X_{t}-x)b(X_{t})\mathrm{d}t}{\frac{1}{T}\int_{0}^{T}K_{h}(X_{t}-x)\mathrm{d}t}+\frac{\frac{1}{T}\int_{0}^{T}K_{h}(X_{t}-x)\sigma\mathrm{d}W_{t}}{\frac{1}{T}\int_{0}^{T}K_{h}(X_{t}-x)\mathrm{d}t}\\
=&\frac{G_{T,1}(x)+G_{T,2}(x)}{F_{T}(x)}.
\end{aligned}
\end{equation*}

For the convenience of the following discussion, we do some marks. Let 
\begin{equation*}
\begin{aligned}
&H_{t,s}=\sigma(W_{t}-W_{s})+(L_{t}-L_{s})-(R_{t}-R_{s}).
\end{aligned}
\end{equation*}

We give the proof of Lemma \ref{lemE} as follows.

\noindent\textbf{Proof of Lemma \ref{lemE}}. a. Note that $b(\cdot)$ is Lipschitz continuous and not identically $0$ on $[l,u]$. It suffices to verify that $X$ satisfies Conditions $5.2-5.5$ in \cite{Budhiraja2007}. Hence, $X$ has a unique invariant probability measure.  Define $m(x)$ and $s(x)$ the scale and speed densities
\begin{equation*}
	s(x)=\exp \left(\frac{2}{\sigma^{2}} \int_{l}^{x} b(r) \mathrm{d} r\right), \quad \text { and } \quad m(x)=\frac{2}{\sigma^{2} s(x)}.
\end{equation*}
If there exists a stationary density $\psi(x)$, it satisfies the Kolmogorov backward equation (\cite{Han2016,Karlin1981})
\begin{equation*}
	\frac{\sigma^{2}}{2} \frac{\mathrm{d}^{2}}{\mathrm{d} x^{2}}(\psi(x))+\frac{\mathrm{d}}{\mathrm{d} x}(b(x) \psi(x))=0 .
\end{equation*}
Solving this equation yeilds 
\begin{equation*}
	\psi(x) =m(x)\left(C_{1} \int_{l}^{x} s(y) \mathrm{d} y+C_{2}\right),
\end{equation*}
where we choose $C_{1}=0$ and $C_{2}=\big(\int_{l}^{u}m(x)\mathrm{d}x\big)^{-1}$. Then 
\begin{equation*}
	\pi_{l,u}(x)=\frac{e^{-\frac{2}{\sigma^{2}} \int_{l}^{x} b(y) \mathrm{d} y}}{\int_{l}^{u} e^{-\frac{2}{\sigma^{2}} \int_{l}^{x} b(y) \mathrm{d} y} \mathrm{d} x}.
\end{equation*}

\noindent b. From that $X$ has a unique invariant probability measure, we can conclude that $\{X_{t},t\geq0\}$ is ergodic (\cite{Han2016}).

\noindent c. The $\Delta$-skeleton $\{X_{t_{k}},k=0,1,\cdots,n\}$ is ergodic follows from the Theorem $16.0.1$ in \cite{Meyn1993}.
\hfill$\square$

Now, we prepare some preliminary lemmas.

\begin{lem}\label{lemF}
Under conditions \ref{asp1}-\ref{aspl}, we have
\begin{equation*}
F_{n}(x)\rightarrow F(x) \quad \text{and} \quad F_{T}(x)\rightarrow F(x),
\end{equation*}
as $n, T\rightarrow\infty$.
\end{lem}
\textbf{Proof of Lemma \ref{lemF}}. Note that $K_{h}\in L_{1}([l,u])$. Then the desired results are obatained immediately by Lemma \ref{lemE}.
\hfill$\square$

\begin{lem}\label{lemD}
Under conditions \ref{asp1}-\ref{aspl}, we have
\begin{equation*}
\sup_{t>s\in\mathbb{R}^{+}}|X_{t}-X_{s}|\leq \big(b(X_{s})(t-s)+|H_{t,s}|\big)e^{M(t-s)},
\end{equation*}
where $M$ is the Lipschitz constant.
\end{lem}
\noindent\textbf{Proof of Lemma \ref{lemD}}. Note that
\begin{equation*}
X_{t}-X_{s}=\int_{s}^{t}b(X_{u})\mathrm{d}u+H_{t,s}.
\end{equation*}
By the Lipschitz continuity of $b$, we have
\begin{equation*}
\begin{aligned}
|X_{t}-X_{s}|&\leq \int_{s}^{t}b(X_{u})\mathrm{d}u+|H_{t,s}|\\
&\leq \int_{s}^{t}\big(|b(X_{u})-b(X_{s})|+b(X_{u})\big)\mathrm{d}u+|H_{t,x}|.
\end{aligned}
\end{equation*}
By the Lipschitz continuity of $b$ again, we have
\begin{equation*}
|X_{t}-X_{s}|\leq \int_{s}^{t}M|X_{u}-X_{s}|\mathrm{d}u+b(X_{s})(t-s)+|H_{t,s}|.
\end{equation*}
From Gronwall's inequality, we obtain the desired results.
\hfill$\square$

\begin{lem}\label{lemWC}
Under conditions \ref{asp1}-\ref{aspl}, we have
\begin{equation*}
G_{T,2}(x)\stackrel{P}{\longrightarrow}0,
\end{equation*}
as $T\rightarrow\infty$.
\end{lem}
\noindent\textbf{Proof of Lemma \ref{lemWC}}. Note that 
\begin{equation*}
\mathbb{E}(G_{T,2}(x))=0.
\end{equation*}
By It$\hat{o}$ isometry, we have
\begin{equation*}
\mathbb{E}(G^{2}_{T,2}(x))=\frac{1}{T^{2}}\int_{0}^{T}\mathbb{E}\bigg(K_{h}(X_{t}-x)\bigg)^{2}\mathrm{d}t.
\end{equation*}
Since that $K$ is bounded and integrable on $[l,u]$, by the dominated convergence theorem, we have
\begin{equation*}
\lim_{T\rightarrow\infty}\mathbb{E}(G^{2}_{T,2}(x))=\lim_{T\rightarrow\infty}\frac{1}{T^{2}}\int_{0}^{T}\bigg(K_{h}(X_{t}-x)\bigg)^{2}\mathrm{d}t.
\end{equation*}
By Lemma \ref{lemE}, we have
\begin{equation*}
\begin{aligned}
\lim_{T\rightarrow\infty}\mathbb{E}(G^{2}_{T,2})=&\frac{1}{T}\int_{l}^{u}\big(K_{h}(y-x)\big)^{2}\pi_{l,u}(y)\mathrm{d}y\\
=&\frac{1}{Th}\int_{-1}^{1}K^{2}(r)\pi_{l,u}(x+hr)\mathrm{d}r\\
=&O\big((Th)^{-1}\big).
\end{aligned}
\end{equation*}
From Chebyshev's inequality and for $\epsilon\in(-\frac{1}{2},0)$, we have
\begin{equation*}
\begin{aligned}
\mathbb{P}(\lim_{T\rightarrow\infty}G_{T,2}(x)>(Th)^{\epsilon})\leq \frac{O\big((Th)^{-1}\big)}{(Th)^{2\epsilon}}=O\big((Th)^{-1-2\epsilon}\big),
\end{aligned}
\end{equation*}
which goes to $0$ as $T\rightarrow\infty$.
\hfill$\square$

\begin{lem}\label{lemWD}
Under conditions \ref{asp1}-\ref{aspl}, we have
\begin{equation*}
G_{n,2}(x)\stackrel{P}{\longrightarrow}0,
\end{equation*}
as $n\rightarrow\infty$.
\end{lem}
\noindent\textbf{Proof of Lemma \ref{lemWD}}. Let
\begin{equation}\label{eqgn2}
\phi_{n}(t,x)=\sum_{k=0}^{n-1}K_{h}(X_{t_{k}}-x)\sigma I_{\{t_{k},t_{k+1}\}}(t).
\end{equation}
Hence,
\begin{equation*}
\lim_{n\rightarrow\infty}G_{n,2}(x)=\lim_{n\rightarrow\infty}\frac{1}{n\Delta}\int_{0}^{t_{n}}\phi_{n}(t,x)\mathrm{d}W_{t}.
\end{equation*}
By some similar arguments as in the proof of Lemma \ref{lemWC}, we have
\begin{equation*}
\lim_{n\rightarrow\infty}\frac{1}{n\Delta}\int_{0}^{t_{n}}\phi_{n}(t,x)\mathrm{d}W_{t}\stackrel{P}{\longrightarrow} 0,
\end{equation*}
which completes the proof.
\hfill$\square$

\noindent\textbf{Proof of Theorem \ref{thmConsD}}. Note that
\begin{equation*}
\begin{aligned}
\hat{b}_{n}(x)-b(x)=&\frac{\frac{1}{n\Delta}\sum_{k=0}^{n-1}K_{h}(X_{t_{k}}-x)\int_{t_{k}}^{t_{k+1}}\big(b(X_{t_{k}})-b(X_{t_{k}})\big)\mathrm{d}t}{\frac{1}{n}\sum_{k=0}^{n-1}K_{h}(X_{t_{k}}-x)}\\
&+\frac{\frac{1}{n}\sum_{k=0}^{n-1}K_{h}(X_{t_{k}}-x)\big(b(X_{t_{k}})-b(x)\big)}{\frac{1}{n}\sum_{k=0}^{n-1}K_{h}(X_{t_{k}}-x)}\\
&+\frac{\frac{\sigma}{n\Delta}\sum_{k=0}^{n-1}K_{h}(X_{t_{k}}-x)\delta W_{t_{k}}}{\frac{1}{n}\sum_{k=0}^{n-1}K_{h}(X_{t_{k}}-x)}\\
=&\frac{G_{n,3}(x)+G_{n,4}(x)+G_{n,2}(x)}{F_{n}(x)}.
\end{aligned}
\end{equation*}

By Lemma \ref{lemF} and \ref{lemWD}, we have
\begin{equation*}
\frac{G_{n,2}}{F_{n}(x)}\stackrel{P}{\longrightarrow}0,
\end{equation*}
as $n\rightarrow\infty$. 

By the Lipschitz continuous of $b$ and Lemma \ref{lemD}, we have
\begin{equation*}
\begin{aligned}	G_{n,3}(x)\leq&\frac{1}{nh\Delta}\sum_{k=0}^{n-1}K(\frac{X_{t_{k}}-x}{h})\int_{t_{k}}^{t_{k+1}}M(X_{t}-X_{t_{k}})\mathrm{d}t\\
\leq&\frac{1}{nh}\sum_{k=0}^{n-1}K(\frac{X_{t_{k}}-x}{h})\big(b(X_{t_{k}})\Delta+|H_{t_{k+1},t_{k}}|\big)e^{M\Delta}.
\end{aligned}
\end{equation*}
By the properties of the two processes $L$ and $R$, we have
\begin{equation*}
\begin{aligned}
&\delta L_{t_{k}}=\max\big(0,A_{t_{k}}-(X_{t_{k}}-l)\big),\\
&\delta R_{t_{k}}=\max\big(0,B_{t_{k}}+(X_{t_{k}}-u)\big),
\end{aligned}
\end{equation*}
where
\begin{equation*}
\begin{aligned}
&A_{t_{k}}=\sup_{s\in[t_{k},t_{k+1})}\big\{-b(X_{t_{k}})(s-t_{k})-\sigma(W_{s}-W_{t_{k}})\big\}\\
&B_{t_{k}}=\sup_{s\in[t_{k},t_{k+1})}\big\{b(X_{t_{k}})(s-t_{k})+\sigma(W_{s}-W_{t_{k}})\big\}.
\end{aligned}
\end{equation*}
By the local H\"older continuity of Brownian motion, we have
\begin{equation*}
H_{t_{k+1},t_{k}}\leq C_{1}\Delta^{\frac{1}{2}-\epsilon},
\end{equation*}
where $\epsilon$ is a small enough positive number and $C_{1}$ are some constants. Hence 
\begin{equation}\label{eqgn3}
\begin{aligned}
G_{n,3}(x)&\leq \frac{1}{nh}\sum_{k=0}^{n-1}K(\frac{X_{t_{k}}-x}{h})\big(b(X_{t_{k}})\Delta+\sup_{|t-s|\leq\Delta}|H_{t,s}|\big)e^{M\Delta}\\
&\leq \frac{1}{h}\|K\|_{\infty}(\|b\|_{\infty}+C_{1})\Delta^{\frac{1}{2}-\epsilon}\\
&=O(h^{-1}\Delta^{\frac{1}{2}-\epsilon}),
\end{aligned}
\end{equation}
which goes to $0$, as $n\rightarrow\infty$.
Hence, we have $G_{n,3}(x)/F_{n}(x)\rightarrow 0$ almost surely, as $n\rightarrow\infty$.

By the properties of the kernel function $K$, we have
\begin{equation*}
K(\frac{y-x}{h})=K(\frac{y-x}{h})I_{\{[x-h,x+h]\}}(y).
\end{equation*}
Let $Z_{1}=\{k, X_{t_{k}}\in[x-h,x+h]\}$ and $Z_{2}=\{k, X_{t_{k}}\notin[x_h,x+h]\}$.
Then,
\begin{equation*}
\begin{aligned}
G_{n,4}(x)=&\frac{1}{n}\sum_{Z_{1}}K_{h}(X_{t_{k}}-x)\big(b(X_{t_{k}})-b(x)\big)\\
&+\frac{1}{n}\sum_{Z_{2}}K_{h}(X_{t_{k}}-x)\big(b(X_{t_{k}})-b(x)\big),
\end{aligned}
\end{equation*}
where the second term is $0$. By the Lipschitz continuous of $b$, we have
\begin{equation}\label{eqgn4}
\begin{aligned}
&\frac{1}{n}\sum_{Z_{1}}K_{h}(X_{t_{k}}-x)\big(b(X_{t_{k}})-b(x)\big)\\
\leq&\frac{1}{n}\sum_{Z_{1}}K_{h}(X_{t_{k}}-x)M(X_{t_{k}}-x)\\
\leq&F_{n}(x)Mh,
\end{aligned}
\end{equation}
which goes to $0$, as $h\rightarrow0$ and $n\rightarrow\infty$. Thus, $G_{n,4}(x)/F_{n}(x)\rightarrow 0$ almost surely, as $n\rightarrow\infty$.
\hfill$\square$

\noindent\textbf{Proof of Theorem \ref{thmConsC}}. Note that
\begin{equation*}
\begin{aligned}
\tilde{b}_{T}(x)-b(x)=&\frac{\frac{1}{T}\int_{0}^{T}K_{h}(X_{t}-x)\big(b(X_{t})-b(x)\big)\mathrm{d}t}{\frac{1}{T}\int_{0}^{T}K_{h}(X_{t}-x)\mathrm{d}t}+\frac{\frac{1}{T}\int_{0}^{T}K_{h}(X_{t}-x)\sigma\mathrm{d}W_{t}}{\frac{1}{T}\int_{0}^{T}K_{h}(X_{t}-x)\mathrm{d}t}\\
=&\frac{G_{T,3}(x)+G_{T,2}(x)}{F_{T}(x)}.
\end{aligned}
\end{equation*}
By Lemma \ref{lemF} and \ref{lemWC}, we have
\begin{equation*}
\frac{G_{T,2}(x)}{F_{T}(x)}\stackrel{P}{\longrightarrow}0,
\end{equation*}
as $T\rightarrow\infty$.

Let $Z_{1}^{\prime}=\{t,X_{t}\in[x-h,x+h]\}$ and $Z_{2}^{\prime}=\{t, X_{t}\notin[x-h,x+h]\}$. Then,
\begin{equation*}
\begin{aligned}
G_{T,3}(x)=&\frac{1}{T}\int_{0}^{T}K_{h}(X_{t}-x)\big(b(X_{t})-b(x)\big)\mathrm{d}t\\
=&\frac{1}{T}\int_{Z_{1}^{\prime}}K_{h}(X_{t}-x)\big(b(X_{t})-b(x)\big)\mathrm{d}t\\&+\frac{1}{T}\int_{Z_{2}^{\prime}}K_{h}(X_{t}-x)\big(b(X_{t})-b(x)\big)\mathrm{d}t,
\end{aligned}
\end{equation*}
where the second term is $0$. By the Lipschitz continuous of $b$, we have
\begin{equation}\label{eqgt3}
\begin{aligned}
&\int_{Z_{1}^{\prime}}K_{h}(X_{t}-x)\big(b(X_{t})-b(x)\big)\mathrm{d}t\\
\leq&\int_{Z_{1}^{\prime}}K_{h}(X_{t}-x)Mh\mathrm{d}t
\\
\leq&F_{T}(x)Mh.
\end{aligned}
\end{equation}
Thus, $G_{T,3}(x)/F_{T}(x)$ goes to $0$ as $h\rightarrow0$ and $n\rightarrow\infty$.
\hfill$\square$

\noindent\textbf{Proof of Theorem \ref{thmAD}}. Note that
\begin{equation*}
\begin{aligned}
\sqrt{nh\Delta}\big(\hat{b}_{n}(x)-b(x)\big)=&\frac{\sqrt{nh\Delta}\big(G_{n,3}(x)+G_{n,4}(x)+G_{n,2}(x)\big)}{F_{n}(x)}.
\end{aligned}
\end{equation*}
By Eq. (\ref{eqgn3}), we have 
\begin{equation*}
\sqrt{nh\Delta}G_{n,3}(x)=\sqrt{nh\Delta}O(h^{-1}\Delta^{\frac{1}{2}-\epsilon})=O(\sqrt{nh^{-1}\Delta^{2-2\epsilon}}),
\end{equation*}
which goes to $0$ as $n\rightarrow\infty$.
By Eq. (\ref{eqgn4}), we have
\begin{equation*}
\sqrt{nh\Delta}G_{n,3}(x)=\sqrt{nh\Delta}O(h)=O(\sqrt{nh^{3}\Delta}),
\end{equation*}
which goes to $0$ as $n\rightarrow\infty$. 
Then, it suffices to show that $\sqrt{nh\Delta}G_{n,2}(x)$ comvergences in law to a centered normal distribution, as $n\rightarrow\infty$.
From Eq. \ref{eqgn2}, we have
\begin{equation*}
\begin{aligned}
\sqrt{nh\Delta}G_{n,2}(x)=\sqrt{\frac{h}{n\Delta}}\int_{0}^{t_{n}}\phi_{n}(t,x)\mathrm{d}W_{t},
\end{aligned}
\end{equation*}
which is a Gaussian random variable.
Moreover, $\mathbb{E}\big(\sqrt{nh\Delta}G_{n,2}(x)\big)=0$.
By It$\hat{o}$ isometry and dominated convergence theorem, we have
\begin{equation*}
\begin{aligned}
\mathbb{E}\big[\big(\sqrt{nh\Delta}G_{n,2}(x)\big)^{2}\big]=&\frac{h}{n\Delta}\int_{0}^{t_{n}}\mathbb{E}\big[\big(\phi_{n}(t,x)\big)^{2}\big]\mathrm{d}t\\
=&\frac{h}{n\Delta}\int_{0}^{t_{n}}\big(\phi_{n}(t,x)\big)^{2}\mathrm{d}t\\
=&\frac{h\sigma^{2}}{n}\sum_{k=0}^{n-1}\big(K_{h}(X_{t_{k}}-x)\big)^{2}.
\end{aligned}
\end{equation*}
By Lemma \ref{lemE}, we have
\begin{equation*}
\begin{aligned}
	\lim_{n\rightarrow\infty}\mathbb{E}\big[\big(\sqrt{nh\Delta}G_{n,2}(x)\big)^{2}\big]=&\lim_{n\rightarrow\infty}\frac{h\sigma^{2}}{n}\sum_{k=0}^{n-1}\big(K_{h}(X_{t_{k}}-x)\big)^{2}\\
	=&h\sigma^{2}\int_{l}^{u}\big(K_{h}(y-x)\big)^{2}\pi_{l,u}(y)\mathrm{d}y\\
	=&\sigma^{2}\int_{\frac{u-x}{h}}^{\frac{l-x}{h}}K^{2}(u)\pi_{l,u}(x+hr)\mathrm{d}r\\
	=&\sigma^{2}F(x),
\end{aligned}
\end{equation*}
which implies that $\sqrt{nh\Delta}G_{n,2}(x)\sim\mathcal{N}\big(0,\sigma^{2}F(x)\big)$. By Lemma \ref{lemF} and Slutsky' theorem, we have
\begin{equation*}
\frac{\sqrt{nh\Delta}G_{n,2}(x)}{F_{n}(x)}\sim\mathcal{N}\big(0,\sigma^{2}F^{-1}(x)\big),
\end{equation*}
which completes the proof.
\hfill$\square$

\noindent\textbf{Proof of Theorem \ref{thmAC}}. Note that
\begin{equation*}
\sqrt{Th}\big(\tilde{b}_{T}(x)-b(x)\big)=\frac{\sqrt{Th}\big(G_{T,3}+G_{T,2}\big)}{F_{T}(x)}.
\end{equation*}
By Eq. (\ref{eqgt3}), we have
\begin{equation*}
\sqrt{Th}G_{T,3}(x)=O(\sqrt{Th^{3}}),
\end{equation*}
which goes to $0$ as $T\rightarrow\infty$.
A similar proof of the proof of Theorem \ref{thmAD}, we have
\begin{equation*}
\sqrt{Th}G_{T,2}(x)\sim\mathcal{N}\big(0, \sigma^{2}F(x)\big).
\end{equation*}
By Lemma \ref{lemF} and Slutsky' theorem, we have
\begin{equation*}
	\frac{\sqrt{Th}G_{n,2}(x)}{F_{T}(x)}\sim\mathcal{N}\big(0,\sigma^{2}F^{-1}(x)\big),
\end{equation*}
which completes the proof.
\hfill$\square$

\section{Several Examples and Numerical Results}\label{sec4}
\subsection{Several examples}
In this subsection, we present several examples on the cases of two-sided barriers and one-sided barrier for illustration. 
The diffusion parameter $\sigma=0.2$ and the drift function we choose shall satisfy the condition \ref{asp1} and \ref{asp2}. Hence, we considered the following cases 
\begin{enumerate}[1.]
	\item $b_{1}(x)=\sin(2\pi x)+1.5x$
	\item $b_{2}(x)=\sqrt{1+x^{2}}$
	\item $b_{3}(x)=2\sqrt{x}$.
\end{enumerate}
We plug the different cases of $b$ into model \ref{eqmodelT} and \ref{eqmodelo} to test the proposed estimators. Moreover, we take the lower reflecting barrier $l=0$ and the upper reflecting barrier $u=3$.

\subsection{Numerical results}
In this subsection, we examine the numerical behavior of the proposed estimator. For a simulation of reflected stochastic processes, we make use of the numerical method presented in \cite{Lepingle1995}, which is known to yield the same rate of convergence as the usual Euler-Maruyama method. Note that there is always the discretely observed process. We only do simulations for discrete-type estimator. 

For each setting, we generate $N=1000$ Monte Carlo simulations of the sample paths, each consisting of $n=400, 900$ or $1600$ observations. Based on conditions \ref{asp4} and \ref{aspamd}, to obtain a N-W estimator with asymptotically negligible bias, we employ $\Delta=n^{-2/3}$ and bandwidths in the range $(n^{-0.3},n^{-0.15})$. The kernel function is the Epanechnikov kernel, which is $K(x)=0.75(1-x^{2})_{+}$. Further simulations evidence that the use of other kernel function has little impact on the estimator's empirical performance, we ommit the results here. 

Let $\{x_{i}, i=1,\cdots,300\}$ be the estimated set where the $300$ points in the set are taken uniformly from $[0,3]$.
The proposed N-W estimator is evaluated by the square-Root of Average Square Errors (RASE), the sample standard deviation (Std.dev) and median (Median) of RASE, where
\begin{equation*}
RASE=\bigg[\frac{1}{300}\sum_{i=1}^{300}\big(\hat{b}(x_{i})-b(x_{i})\big)^{2}\bigg]^{\frac{1}{2}}.
\end{equation*} 

\begin{figure}
\centering
\subfigure[$b(x)=\sin(2\pi x)+1.5x$ with two-sided reflecting barrier.]{
    \includegraphics[scale=0.5]{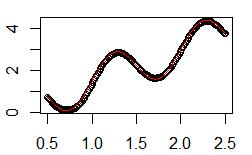}
}
\subfigure[$b(x)=\sin(2\pi x)+1.5x$ with one-sided reflecting barrier.]{
    \includegraphics[scale=0.5]{two1600.3b1}
}
\subfigure[$b(x)=\sqrt{1+x^{2}}$ with two-sided reflecting barrier.]{
	\includegraphics[scale=0.5]{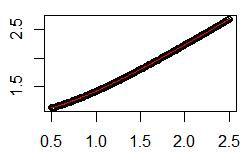}
}
\subfigure[$b(x)=\sqrt{1+x^{2}}$ with one-sided reflecting barrier.]{
	\includegraphics[scale=0.5]{two1600.3b2}
}
\subfigure[$b(x)=2\sqrt{x}$ with two-sided reflecting barrier.]{
	\includegraphics[scale=0.5]{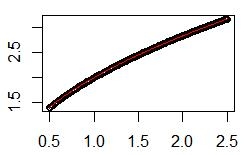}
}
\subfigure[$b(x)=2\sqrt{x}$ with one-sided reflecting barrier.]{
	\includegraphics[scale=0.5]{two1600.3b3}
}
\caption{The N-W estimator of $b(\cdot)$ with $n=1600$ and $h=n^{-0.3}$, where the dot is the estimator and the solid red line is the true value.}
\label{fig1}
\end{figure}

Table \ref{table1}-\ref{table3} summarize the main findings over 1000 simulations. We observe that as the sample size increases or the bandwidth decreases, the RASE decreases and is small, that the empirical and model-based standard errors agree reasonably. The performance improves with larger sample sizes and smaller bandwidths. Figure \ref{fig1} represents the N-W estimator $\hat{b}(\cdot)$ from the three cases with $n=1600$ and $h=n^{-0.3}$. It shows that the N-W estimator performs reasonably well.

\begin{table}
	\caption{Simulation results with different sample size $n$ and bandwidth $h$.}
	\begin{tabular}{cccccccc}
		\hline 
		\multicolumn{8}{c}{Case 1: $b_{1}(x)=\sin(2\pi x)+1.5x$.} \tabularnewline
		\hline 
		\hline 
	    &&\multicolumn{3}{c}{Two-sided}
	    &\multicolumn{3}{c}{One-sided} 
	    \tabularnewline
        \cline{3-5}  \cline{6-8}
		\makebox[0.06\textwidth]{n}&
		\makebox[0.06\textwidth]{BD} &
		\makebox[0.1\textwidth]{RASE}&  
		\makebox[0.1\textwidth]{Std.dev}&
		\makebox[0.1\textwidth]{Median}&
		\makebox[0.1\textwidth]{RASE}&  
		\makebox[0.1\textwidth]{Std.dev}&
		\makebox[0.1\textwidth]{Median}
		\tabularnewline
		400  & $n^{-0.3}$  & 0.190 & 0.003 & 0.205 
		                   & 0.128 & 0.004 & 0.150\tabularnewline
		     & $n^{-0.2}$  & 0.438 & 0.003 & 0.443 
		                   & 0.338 & 0.004 & 0.345\tabularnewline
		     & $n^{-0.15}$ & 0.619 & 0.003 & 0.622
		                   & 0.515 & 0.004 & 0.519\tabularnewline
		900  & $n^{-0.3}$  & 0.135 & 0.003 & 0.150 
		                   & 0.085 & 0.003 & 0.109\tabularnewline
		     & $n^{-0.2}$  & 0.365 & 0.003 & 0.368 
		                   & 0.264 & 0.003 & 0.269\tabularnewline
		     & $n^{-0.15}$ & 0.553 & 0.003 & 0.555
		                   & 0.438 & 0.003 & 0.442\tabularnewline
		1600 & $n^{-0.3}$  & 0.105 & 0.002 & 0.120
		                   & 0.064 & 0.003 & 0.088\tabularnewline
		     & $n^{-0.2}$  & 0.315 & 0.002 & 0.318
		                   & 0.220 & 0.003 & 0.225\tabularnewline
		     & $n^{-0.15}$ & 0.505 & 0.002 & 0.507 
		                   & 0.389 & 0.003 & 0.392\tabularnewline
		\hline 
	\end{tabular}
	\label{table1}
\end{table}

\begin{table}
	\caption{Simulation results with different sample size $n$ and bandwidth $h$.}
	\begin{tabular}{cccccccc}
		\hline 
		\multicolumn{8}{c}{Case 2: $b_{2}(x)=\sqrt{1+x^{2}}$.} \tabularnewline
		\hline 
		\hline 
		&&\multicolumn{3}{c}{Two-sided}
		&\multicolumn{3}{c}{One-sided} 
		\tabularnewline
		\cline{3-5}  \cline{6-8}
		\makebox[0.06\textwidth]{n}&
		\makebox[0.06\textwidth]{BD} &
		\makebox[0.1\textwidth]{RASE}&  
		\makebox[0.1\textwidth]{Std.dev}&
		\makebox[0.1\textwidth]{Median}&
		\makebox[0.1\textwidth]{RASE}&  
		\makebox[0.1\textwidth]{Std.dev}&
		\makebox[0.1\textwidth]{Median}
		\tabularnewline
		400  & $n^{-0.3}$  & 0.017 & 0.003 & 0.070 
		                   & 0.003 & 0.003 & 0.071\tabularnewline
             & $n^{-0.2}$  & 0.043 & 0.002 & 0.064
                           & 0.006 & 0.002 & 0.530\tabularnewline
             & $n^{-0.15}$ & 0.067 & 0.001 & 0.078 
                           & 0.010 & 0.002 & 0.047\tabularnewline
        900  & $n^{-0.3}$  & 0.013 & 0.003 & 0.060
                           & 0.003 & 0.003 & 0.061\tabularnewline
             & $n^{-0.2}$  & 0.037 & 0.001 & 0.055 
                           & 0.004 & 0.002 & 0.044\tabularnewline
             & $n^{-0.15}$ & 0.060 & 0.001 & 0.069 
                           & 0.008 & 0.002 & 0.038\tabularnewline
        1600 & $n^{-0.3}$  & 0.011 & 0.002 & 0.054
                           & 0.003 & 0.002 & 0.055\tabularnewline
             & $n^{-0.2}$  & 0.032 & 0.001 & 0.048
                           & 0.003 & 0.002 & 0.038\tabularnewline
             & $n^{-0.15}$ & 0.055 & 0.001 & 0.062 
                           & 0.007 & 0.001 & 0.033\tabularnewline
		\hline 
	\end{tabular}
	\label{table2}
\end{table}

\begin{table}
	\caption{Simulation results with different sample size $n$ and bandwidth $h$.}
	\begin{tabular}{cccccccc}
		\hline 
		\multicolumn{8}{c}{Case 3: $b_{3}(x)=2\sqrt{x}$.} \tabularnewline
		\hline 
		\hline 
		&&\multicolumn{3}{c}{Two-sided}
		&\multicolumn{3}{c}{One-sided} 
		\tabularnewline
		\cline{3-5}  \cline{6-8}
		\makebox[0.06\textwidth]{n}&
		\makebox[0.06\textwidth]{BD} &
		\makebox[0.1\textwidth]{RASE}&  
		\makebox[0.1\textwidth]{Std.dev}&
		\makebox[0.1\textwidth]{Median}&
		\makebox[0.1\textwidth]{RASE}&  
		\makebox[0.1\textwidth]{Std.dev}&
		\makebox[0.1\textwidth]{Median}
		\tabularnewline
		400  & $n^{-0.3}$  & 0.036 & 0.003 & 0.084
		                   & 0.035 & 0.003 & 0.086\tabularnewline
		     & $n^{-0.2}$  & 0.065 & 0.003 & 0.086
		                   & 0.060 & 0.003 & 0.085\tabularnewline
		     & $n^{-0.15}$ & 0.089 & 0.003 & 0.101
		                   & 0.079 & 0.003 & 0.096\tabularnewline
		900  & $n^{-0.3}$  & 0.029 & 0.003 & 0.072
		                   & 0.028 & 0.003 & 0.073\tabularnewline
		     & $n^{-0.2}$  & 0.055 & 0.003 & 0.073 
		                   & 0.051 & 0.003 & 0.071\tabularnewline
		     & $n^{-0.15}$ & 0.079 & 0.002 & 0.088
		                   & 0.069 & 0.003 & 0.082\tabularnewline
		1600 & $n^{-0.3}$  & 0.025 & 0.003 & 0.064 
		                   & 0.024 & 0.003 & 0.066\tabularnewline
		     & $n^{-0.2}$  & 0.049 & 0.002 & 0.064
		                   & 0.046 & 0.002 & 0.063\tabularnewline
		     & $n^{-0.15}$ & 0.072 & 0.002 & 0.080
		                   & 0.064 & 0.002 & 0.074\tabularnewline
		\hline  
	\end{tabular}
	\label{table3}
\end{table}

\section{Conclusion}\label{sec5}
In this papaer, we study the N-W estimators for the drift function of two-sided reflected processes. We proposed the discrete-type and continuous-type N-W estimators based on discretely observed processes and continuously observed processes respectively. 
We obatin the consistency of the proposed estimators and establish their asymptotic distributions. Our method can be applied to the reflected stochastic processes with only one-sided reflecting barrier. A simulation study concerning the finite sample behavior has also been presented. The numerical results show that the proposed estimators are adequade for practical use.

\end{document}